\documentclass[reqno]{amsart}
\usepackage[a4paper]{geometry} 
\usepackage{amsmath}                                                            
\usepackage{amsfonts}
\usepackage{amsthm}
\usepackage{listings}
\usepackage{xcolor}
\definecolor{dkgreen}{rgb}{0,0.6,0}
\definecolor{mauve}{rgb}{0.58,0,0.82}
\usepackage{url}
\usepackage{caption}
\usepackage{listings}
\usepackage{color}
\usepackage{inconsolata}
\usepackage{a4wide}
\usepackage{bm}
\usepackage{appendix}
\usepackage{tikz}
\usetikzlibrary{calc}

\usepackage[labelformat=simple]{subcaption}

\usepackage{vem_macros}

\usepackage{algorithm}
\usepackage[noend]{algpseudocode}
\numberwithin{equation}{section}

\makeatletter
\def\BState{\State\hskip-\ALG@thistlm}
\makeatother

\definecolor{mygreen}{rgb}{0,0.6,0}
\definecolor{mygray}{rgb}{0.5,0.5,0.5}
\definecolor{mymauve}{rgb}{0.58,0,0.82}
\definecolor{cream}{rgb}{1.0, 1.0, 0.94}

\renewcommand{\PE}[1]{\mathcal{P}_{\E}}
\newcommand{\Pe}[1]{\mathcal{P}_{\edge}}
\renewcommand{\ME}[1]{\mathcal{M}_{\E}}
\renewcommand{\Po}[1]{\Pi^{\E}_0}
\renewcommand{\l}{\ell}
\newcommand{\lh}{\ell_{h}}
\newcommand{\bi}{\varphi_i}
\newcommand{\bj}{\varphi_j}
\newcommand{\biE}{\bi^{\E}}
\newcommand{\bjE}{\bj^{\E}}
\newcommand{\stiffmat}{K}
\newcommand{\stiffmatE}{\stiffmat^{\E}}

\newcommand{\fvec}{F}
\newcommand{\fvecE}{\fvec^{\E}}
\newcommand{\U}{U}
\renewcommand{\Pn}{\Pi^{\E}}
\newcommand{\dimPE}{N_{\mathcal{P}}}
\newcommand{\Dmat}{D}
\newcommand{\projlhsmat}{G}
\newcommand{\projmat}{\bm{\Pi}}
\newcommand{\projrhs}{B}
\newcommand{\mb}{m_{\beta}}

\newcommand{\codelocation}{\url{www2.le.ac.uk/departments/mathematics/research/virtual-element-methods-1/software}}

\title{The virtual element method in 50 lines of Matlab}
\author{Oliver J. Sutton}
\address[O. J. Sutton]{Department of Mathematics, University of Leicester, University Road, Leicester, LE1 7RH, UK}
\email{ojs4@le.ac.uk}

\begin{document}

\numberwithin{lstlisting}{section}

\begin{abstract}
We present a 50-line MATLAB implementation of the lowest order virtual element method for the two-dimensional Poisson problem on general polygonal meshes.
The matrix formulation of the method is discussed, along with the structure of the overall algorithm for computing with a virtual element method.
The purpose of this software is primarily educational, to demonstrate how the key components of the method can be translated into code.
\end{abstract}
\maketitle

\section{Introduction}

The virtual element method, introduced in~\cite{VEIGA:2013wi}, is a generalisation of the standard conforming finite element method for the approximation of solutions to partial differential equations.
The method is designed in such a way as to enable the construction of high order approximation spaces which may include an arbitrary degree of global regularity~\cite{BeiraodaVeiga:2013ii} on meshes consisting of very general polygonal (or polyhedral) elements.
This cocktail of desirable features has attracted the method a lot of attention (see, for example,~\cite{Ahmad:2013uaa,Brezzi:2014dk,deDios:2014uq,Cangiani:2015vh,Cangiani:2016ug,Antonietti:2014kh,Perugia:2015td,Vacca:2015jh}), and is made possible through the \emph{virtual element space} of trial and test functions, which is constructed on each mesh element from functions which are implicitly defined through local PDE problems.
These local problems are designed in such a way that the virtual element space includes a subspace of polynomials of some prescribed degree (referred to as the \emph{degree} of the method) alongside other, typically unknown, \emph{virtual} functions.
In this respect, and like many other approaches to polygonal meshes such as the Polygonal Finite Element Method~\cite{Sukumar:2004dn,Manzini:2014fj} or BEM-based FEM~\cite{Rjasanow:2012bx}, the virtual element method falls within the broad class of Generalised Finite Element Methods~\cite{Strouboulis:2000gs}.
What sets the virtual element method apart from these other approaches, however, is that the extra non-polynomial \emph{virtual} functions never need to be determined or evaluated in practice.
Instead, they are understood and used solely through certain defining properties of the virtual element space and through their degrees of freedom, which, along with the discrete bilinear form, are carefully selected to ensure that the final stiffness matrix can be directly and exactly computed.

There is a history of short, simple, codes being used to demonstrate the practical implementation details of various finite element methods.
We refer, for instance, to the ``Remarks around 50 lines of MATLAB'' paper~\cite{Alberty:1999ff} which presented a simple and transparent MATLAB implementation of the conforming finite element method, the 99-line topology optimisation code presented in~\cite{Sigmund:2001bc}, or the mixed-FEM implementations presented by~\cite{Bahriawati:2005gw}.
Unlike the large workhorse finite element libraries
or commercially available `black-box' finite element software packages, these codes are primarily designed to be educational, demonstrating how the theoretical concepts can be distilled into something practical.
This is a tradition which we continue here, presenting a 50-line MATLAB implementation of the lowest order virtual element method for the Poisson problem in 2D on general polygonal meshes.
To the best of our knowledge, this is the first publicly available implementation of the virtual element method, although there are various references containing details of the matrix formulation of the method, see for instance~\cite{BeiraodaVeiga:2014hz,Cangiani:2015vh,Gain:2013tu}.
In particular,~\cite{BeiraodaVeiga:2014hz} contains detailed explanations of the formulation of the terms in the matrix equations for the high order virtual element method applied to a reaction-diffusion problem.
In many ways, the work here can be seen as a spiritual successor to~\cite{BeiraodaVeiga:2014hz} in the sense that while we restrict ourselves to just the linear virtual element method for the Poisson problem, we take the process one step further and provide a clear, useable MATLAB implementation of the method.

The remainder of this work is structured as follows.
In Section~\ref{sec:modelproblem} we present the model problem of the Poisson problem.
A very brief introduction to the virtual element framework is presented in Section~\ref{sec:method}, with a discussion of the discrete function spaces and bilinear forms.
The details of the implementation of this method are given in Section~\ref{sec:implementation}, where we derive the matrix form of the discrete problem and show how this may be computed in practice.
Section~\ref{sec:usage} contains a brief explanation of how to run the code using MATLAB.
Finally, we offer some concluding remarks and ideas for possible extensions to the code in Section~\ref{sec:conclusion}.
The full code of the implementation is shown in Appendix~\ref{sec:fullcode}, and is 
available to download from:
\begin{center}
\codelocation
\end{center}
alongside several examples of polygonal meshes.

\section{Model problem}
\label{sec:modelproblem}
Let $\Omega \subset \Re^2$ be a polygonal domain and consider the Poisson problem
\begin{equation}
	\begin{split}
		- \Delta u = f \quad &\text{in } \Omega, \\
		u = g \quad &\text{on } \partial \Omega,
	\end{split}
	\label{eq:modelproblem}
\end{equation}
with $f \in L^{2}(\Omega)$ and $g \in H^{1/2}(\partial \Omega)$.
This problem can be written in variational form as: find $u \in H^1_g(\Omega) := \left\{ w \in H^1(\Omega) : w = g \text{ on } \partial \Omega \right\}$ such that
\begin{equation}
	a( u,  v) := (\nabla u, \nabla v) = (f, v) =: \l(v) \qquad \forall \v \in H^1_0(\D)
	\label{eq:variationalproblem}
\end{equation}
where $(\cdot,\cdot)$ denotes the standard $L^2(\D)$ inner product.
This variational problem possesses a unique solution by the Lax-Milgram lemma.

\section{The Virtual Element Method}
\label{sec:method}
Let $\Th$ be a family of partitions of the domain $\D$ into non-overlapping polygonal elements with maximum diameter $h$ whose boundaries are not self-intersecting. 
A theoretical analysis of the method requires certain additional assumptions on the regularity of the mesh, although since we do not present any analysis of the method we do not include these here.
The resulting set of requirements is general enough, however, to include polygonal elements consisting of an arbitrary (but uniformly bounded) number of edges, which may also be non-convex.
However, to simplify the implementation we restrict the mesh to include only elements which contain their own centroid, as defined in~\eqref{eq:centroid}.
We note that this class of elements includes those with co-planar edges, as commonly found in locally refined meshes with hanging nodes, and even non-convex elements.

For a polygonal element $\E$ with $\dimVhE$ edges, we denote its vertices by $\nu_i$ for $i=1,\dots,\dimVhE$, and we adopt the convention that the edge $\e_i$ connects $\nu_i$ and $\nu_{i+1}$, where the indices are understood to wrap within the range 1 to $\dimVhE$.

\subsection{Virtual element function spaces}

The discrete function space is defined to be
\begin{equation*}
	\Vh := \{ \vh \in H^1_g(\D) : \vh|_{\E} \in \VhE \text{ for all } \E \in \Th \}
\end{equation*}
where the local space $\VhE$ on the element $\E$ can be understood through the following three properties:
\begin{itemize}
	\item $\VhE$ includes the space $\PE{1}$ of physical-frame polynomials on $\E$ with total degree $\leq 1$ as a subspace.
	\item Any function in $\VhE$ can be uniquely identified by its values at the vertices of $\E$, which are taken to be the degrees of freedom of the space. 
		We note that this implies that the dimension of the space $\VhE$ is equal to $\dimVhE$.
	\item Every function in $\VhE$ is a linear polynomial on each edge of $\E$.
\end{itemize}
The subspace of linear polynomials provides the approximation power of the virtual element space, and is responsible for the accuracy of the method.
On triangular elements the space consists entirely of these linear polynomials, and thus the method reduces to the standard linear finite element method.
However, on more general shaped polygonal elements the space will also include other, implicitly defined, `virtual' functions, cf.~\eqref{eq:vemspace}.
The method is designed in such a way that these will never need to be explicitly computed or evaluated, and are instead understood solely through their values at the vertices of $\E$, which we take to be the degrees of freedom of the space $\VhE$.
In this respect, the virtual element space can be seen as a straightforward generalisation of the standard linear conforming finite element space on triangles to more general shaped elements.

The first observation we make about this space is that just the properties outlined above allow us to compute the Ritz projection $\Pn : \VhE \rightarrow \PE{1}$ of \emph{any} function in the local virtual element space $\VhE$ onto the subspace of linear polynomials.
This projection is defined for $\vh \in \VhE$ by the conditions
\begin{equation}
	\begin{cases}
		(\nabla (\Pn \vh - \vh), \nabla p)_{0,\E} = 0 \text{ for all } p \in \PE{1}, \\
		\overline{\Pn \vh} = \overline{\vh},
	\end{cases}
	\label{eq:ritzdefn}
\end{equation}
where $\overline{\wh} := \frac{1}{\dimVhE} \sum_{i=1}^{\dimVhE} \wh(\nu_i)$ denotes the average value of $\wh$ at the vertices of $\E$.
This second condition is necessary to fix the constant part of $\Pn\vh$, and is clearly computable for any $\vh \in \VhE$ from just its degrees of freedom.

From~\eqref{eq:ritzdefn}, the divergence theorem, and the fact that the Laplacian of a linear function is zero we have that, for any $\vh \in \VhE$ and $p \in \PE{1}$,
\begin{equation}
	(\nabla \Pn \vh, \nabla p)_{0,\E} = (\nabla \vh, \nabla p)_{0,\E} = \sum_{\edge \in \dE} \inte \vh \n_{\edge} \cdot \nabla p \ds,
	\label{eq:computingPiNabla}
\end{equation}
where $\n_{\edge}$ denotes the unit normal vector to the edge $\edge$ directed out of the element $\E$.
The final expression on the right hand side here can be \emph{exactly} evaluated since $\vh$ is a linear polynomial on each edge of $\E$, entirely determined by its values at the vertices, while the gradient of the linear polynomial $p$ is a known constant.
By picking a basis for the polynomial space $\PE{1}$, equation~\eqref{eq:computingPiNabla} can be written as a matrix problem which can be solved to find the coefficients of $\Pn\vh$ with respect to this polynomial basis.
We will come back to this in Section~\ref{sec:implementation}, although for now we just rely on the fact that this projection is \emph{computable}.

The actual definition of the virtual element space which we use here is the lowest order space introduced in~\cite{VEIGA:2013wi}, given by
\begin{equation}
	\begin{split}
		\VhE := \{ v \in H^1(\E) : \, &\Delta v =0, \quad v|_{\dE} \in C^0(\dE), \\
					 &v|_{\edge} \in \Pe{1} \text{ for each } \edge \in \dE \},
	\end{split}
	\label{eq:vemspace}
\end{equation}
where $\Pe{1}$ denotes the space of linear polynomials on the edge $\edge$.
The fact that the vertex values can be used as degrees of freedom to describe this space is proven in~\cite{VEIGA:2013wi}.

The global degrees of freedom for $\Vh$ are simply taken to be the value of the function at each vertex in the mesh, thus imposing the continuity of the ambient space.
The degrees of freedom at vertices on the domain boundary are fixed in accordance with the boundary condition.
The dimension of the global virtual element space $\Vh$ shall be denoted by $\dimVh$.

\subsection{Discrete bilinear form}
Define the bilinear form $\aE : H^1(\E) \times H^1(\E) \rightarrow \Re$ to be the restriction of $\a$ to the element $\E$, i.e.
	$\aE(\v, \w) := (\nabla \v, \nabla \w)_{0,\E}$
for any $\v, \w \in H^1(\E)$.
Following~\cite{VEIGA:2013wi}, we introduce the discrete counterpart $\ahE : \VhE \times \VhE \rightarrow \Re$ of $\aE$ which we define as
\begin{equation}
	\ahE(\vh, \wh) := (\nabla \Pn \vh, \nabla \Pn \wh)_{0,\E} + \StaE(\vh - \Pn \vh, \wh - \Pn \wh),
	\label{eq:defahE}
\end{equation}
with
\begin{equation*}
	\StaE(\vh, \wh) := \sum_{r=1}^{\dimVhE} \dof_r(\vh) \dof_r(\wh),
\end{equation*}
where $\dof_r(\vh)$ denotes the value of the $r^{\text{th}}$ local degree of freedom defining $\vh$ in $\VhE$ with respect to some arbitrary (but fixed) ordering\footnote{For instance, this could be achieved simply by numbering the vertices of the polygon $\E$.}.
This means that $\StaE$ is simply the Euclidean inner product between vectors of degrees of freedom.
Finally, we define
\begin{equation*}
	\ah(\vh, \wh) := \sum_{\E \in \Th} \ahE(\vh, \wh),
\end{equation*}
to be the discrete counterpart of $\a$.

Crucial to the method is the observation that the local discrete bilinear forms satisfy the following two properties~\cite{VEIGA:2013wi}:
\begin{itemize}
	\item \emph{Polynomial consistency:} for any $\vh \in \VhE$ and $p \in \PE{1}$,
	\begin{equation*}
		\ahE(\vh, p) = \aE(\vh, p).
	\end{equation*}
	\item \emph{Stability:} there exists a constant $\Cstab$ independent of $h$ and $\E$ such that
		\begin{equation*}
			\Cstab^{-1} \aE(\vh, \vh) \leq \ahE(\vh, \vh) \leq \Cstab \aE(\vh, \vh),
		\end{equation*}
		for any $\vh \in \VhE$.
\end{itemize}
The requirement of \emph{polynomial consistency} implies that the method satisfies the patch test commonly used in the engineering literature, expressing the fact that the method is \emph{exact} when the solution is a piecewise linear polynomial with respect to the mesh $\Th$, and provides the accuracy of the method.
The \emph{stability} property, on the other hand, ensures that the discrete bilinear form inherits the continuity and coercivity of the original variational form $\a$, as proven in~\cite{VEIGA:2013wi}.
In the final matrix formulation of the problem, this property can be viewed as ensuring that the problem stiffness matrix is of the correct rank.

In light of these two properties, the two terms of $\ahE$ are referred to as the \emph{consistency} and \emph{stabilising} terms respectively since only the first term is non-zero when either $\vh$ or $\wh$ is a polynomial, and thus single-handedly ensures that the polynomial consistency property is satisfied, while the second term ensures that the stability property is satisfied even when $\vh$ or $\wh$ are in the kernel of $\Pn$.
For a proof that this choice of stabilising term $\StaE$ does indeed satisfy the stability property, we refer to~\cite{Cangiani:2015vh}.

Moreover, both terms of $\ahE$ in~\eqref{eq:defahE} are computable using just the degrees of freedom of the virtual element space (to compute the projector $\Pn$ and to evaluate the stabilising term), and knowledge of the polynomial subspace $\PE{1}$ (to evaluate the consistency term of $\ahE$, which is made of integrals of polynomials, just like in a standard finite element method).
This will be further demonstrated in Section~\ref{sec:implementation}, where it will also become apparent that this particular virtual element method can be implemented \emph{without requiring any quadrature} to compute the stiffness matrix.

Still following~\cite{VEIGA:2013wi}, the linear form $\l$ on the right-hand side of the variational problem~\eqref{eq:variationalproblem} is discretised by $\lh : \VhE \rightarrow \Re$ such that
\begin{equation}
	\lh(\vh) := \sum_{\E \in \Th} (\Po{0} \force, \overline{\vh})_{0,\E},
	\label{eq:deflh}
\end{equation}
where $\Po{0} : \VhE \rightarrow \Re$ denotes the $L^2(\E)$-orthogonal projection onto constants, defined for any $\wh \in \VhE$ to be such that
\begin{equation*}
	\intE (\wh - \Po{0} \wh) \dx = 0.
\end{equation*}

The discrete problem which we solve can then be written as: find $\uh \in \Vh$ such that
\begin{equation}
	\ah(\uh, \vh) = \lh(\vh),
	\label{eq:vemproblem}
\end{equation}
for all $\vh \in \VhE$.

\section{Implementation}
\label{sec:implementation}

As with a typical finite element method, we start by introducing the Lagrangian basis $\{\bi\}_{i=1}^{\dimVh}$ of $\Vh$ with respect to the global set of degrees of freedom, defined by the property that $\bi(\nu_j) = \delta_{ij}$, where $\delta_{ij}$ is the Kronecker delta.
We also introduce the Lagrangian basis of the local virtual element space $\VhE$ on the element $\E$ as $\{\biE\}_{i=1}^{\dimVhE}$, defined by the local equivalent of the same property.

We will also need a basis for the space $\PE{1}$ of local physical frame linear polynomials on each element $\E$.
Many choices are possible here, although in keeping with the convention commonly adopted with virtual element methods, we choose the set of \emph{scaled monomials} of degree 1.
These are defined on the element $\E$ as
\begin{equation}
	\ME{1} := \left\{ \m_1(x,y) := 1, \quad \m_2(x,y) := \frac{x - x_{\E}}{h_{\E}}, \quad \m_3(x,y) := \frac{y - y_{\E}}{h_{\E}} \right\},
	\label{eq:polynomialbasis}
\end{equation}
where $x_{\E}$ and $y_{\E}$ respectively denote the $x$ and $y$ coordinates of the centroid of the element in the standard Cartesian coordinate system, and $h_{\E}$ is the diameter of the element $\E$.
We denote by $\dimPE = 3$ the cardinality of this basis and therefore the dimension of $\PE{1}$

In the hope of avoiding confusion, we adopt the convention of indexing coefficients and basis functions in the basis of $\VhE$ using Latin letters, while those of $\PE{1}$ will be indexed using Greek letters.

With these two bases at our disposal, we can now derive the matrix form of the discrete problem~\eqref{eq:vemproblem}.
Expanding the virtual element solution $\uh$ as
\begin{equation*}
	\uh = \sum_{i = 1}^{\dimVh} U_i \bi,
\end{equation*}
problem~\eqref{eq:vemproblem} can be rewritten using the definitions~\eqref{eq:defahE} of $\ahE$ and~\eqref{eq:deflh} of $\lh$ as: find $U \in \Re^{\dimVh}$ such that
\begin{equation*}
	\sum_{i = 1}^{\dimVh} U_i \sum_{\E \in \Th} \left((\nabla \Pn \bi, \nabla \Pn \bj)_{0,\E} + \StaE(\bi - \Pn \bi, \bj - \Pn \bj) \right) = \sum_{\E} (\Po{0} \force, \overline{\bj})_{0,\E}
\end{equation*}
for $j = 1, \dots, \dimVh$.
This may be expressed in matrix form as
$	\stiffmat \U = \fvec$
where
\begin{equation*}
	\stiffmat_{j,i} = \sum_{\E \in \Th} \left((\nabla \Pn \bi, \nabla \Pn \bj)_{0,\E} + \StaE(\bi - \Pn \bi, \bj - \Pn \bj) \right), \quad \fvec_{j} = \sum_{\E} (\Po{0} \force, \overline{\bj})_{0,\E},
\end{equation*}
for $i,j = 1,\dots,\dimVh$.
Since both these terms are defined through sums over elements, the obvious way to compute the entries of $\stiffmat$ and $\fvec$ is by computing the non-zero local contributions from each element $\E$ in the form of the \emph{local stiffness matrix} $\stiffmatE \in \Re^{\dimVhE \times \dimVhE}$ and \emph{local forcing vector} $\fvecE \in \Re^{\dimVhE}$, given by
\begin{equation}
	\stiffmatE_{j,i} = (\nabla \Pn \biE, \nabla \Pn \bjE)_{0,\E} + \StaE(\biE - \Pn \biE, \bjE - \Pn \bjE), \quad \fvecE_j := (\Po{0} \force, \overline{\bjE})_{0,\E},
	\label{eq:localMatrices}
\end{equation}
for $i,j = 1,\dots,\dimVhE$, and adding them into the corresponding entries of $\stiffmat$ and $\fvec$.

This, of course, dictates that the overall structure of a virtual element method implementation will be much the same as for a standard finite element method, as outlined in Algorithm~\ref{alg:FEM}.
The key point of departure from the standard finite element method is in \emph{how} the element stiffness matrices should be calculated.
Firstly, the computation of the local stiffness matrices relies on first computing the local Ritz projector $\Pn$ on each element.
Secondly, where the implementation of a conventional finite element might rely on a mapping to a reference element, no such equivalent process is possible here because the mesh elements are allowed to be general (possibly non-convex) polygons.
\begin{algorithm}
  \caption{A typical finite/virtual element method implementation}
  \label{alg:FEM}
  \begin{algorithmic}[1]
    \Require Mesh
    \State Initialise stiffness matrix and forcing vector \label{algline:initialisation}
    \For{each element in the mesh} \label{algline:elementloop}
    	\State Compute local stiffness matrix \label{algline:localstiffness}
		\State Compute local forcing contributions \label{algline:localforcing}
		\State Assemble local matrices into the correct places in the global matrices \label{algline:localtoglobal}%
    \EndFor
    \State Condense boundary conditions from matrix system \label{algline:condenseBCs}
    \State Solve the resulting matrix problem \label{algline:solve}
    \State Replace boundary conditions in the solution vector \label{algline:replaceBCs}
    \Ensure Vector containing the degrees of freedom of the discrete solution
  \end{algorithmic}
\end{algorithm}

The full code of the implementation is given in Appendix~\ref{sec:fullcode}.
In the remainder of this section we dissect the code to highlight how the various steps outlined in Algorithm~\ref{alg:FEM} can be expressed in a form which can be implemented in code.
Much of the matrix formulation presented in this section is similar to that in~\cite{BeiraodaVeiga:2014hz}, although here we focus specifically on the case of the lowest order method and take the process a step further to include the details of how each step is accomplished in the code.

\begin{figure}
\centering
\subcaptionbox{Triangles}[4.7cm]{%
\includegraphics[width=4.5cm]{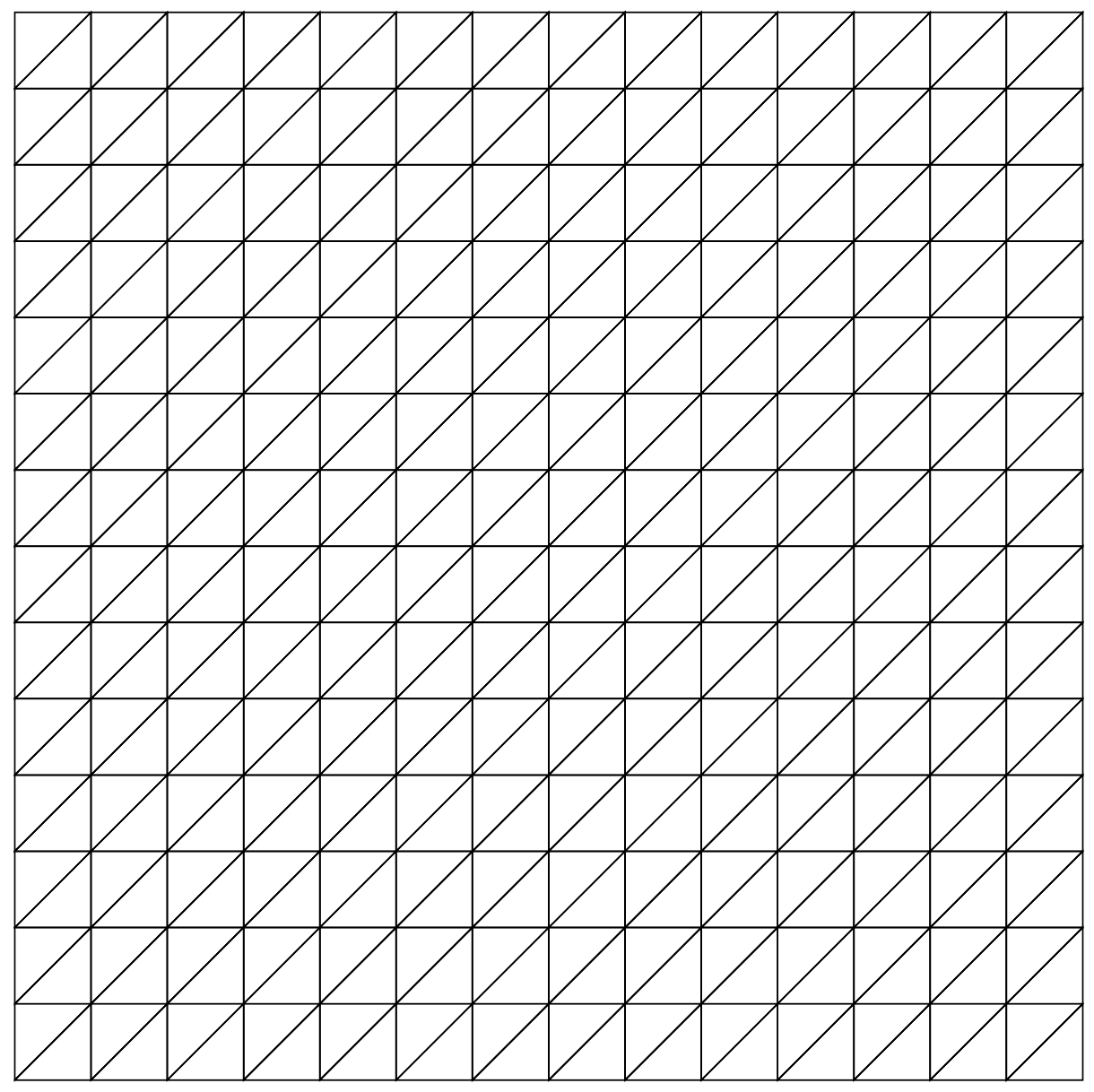}}%
\subcaptionbox{Squares}[4.7cm]{%
\includegraphics[width=4.5cm]{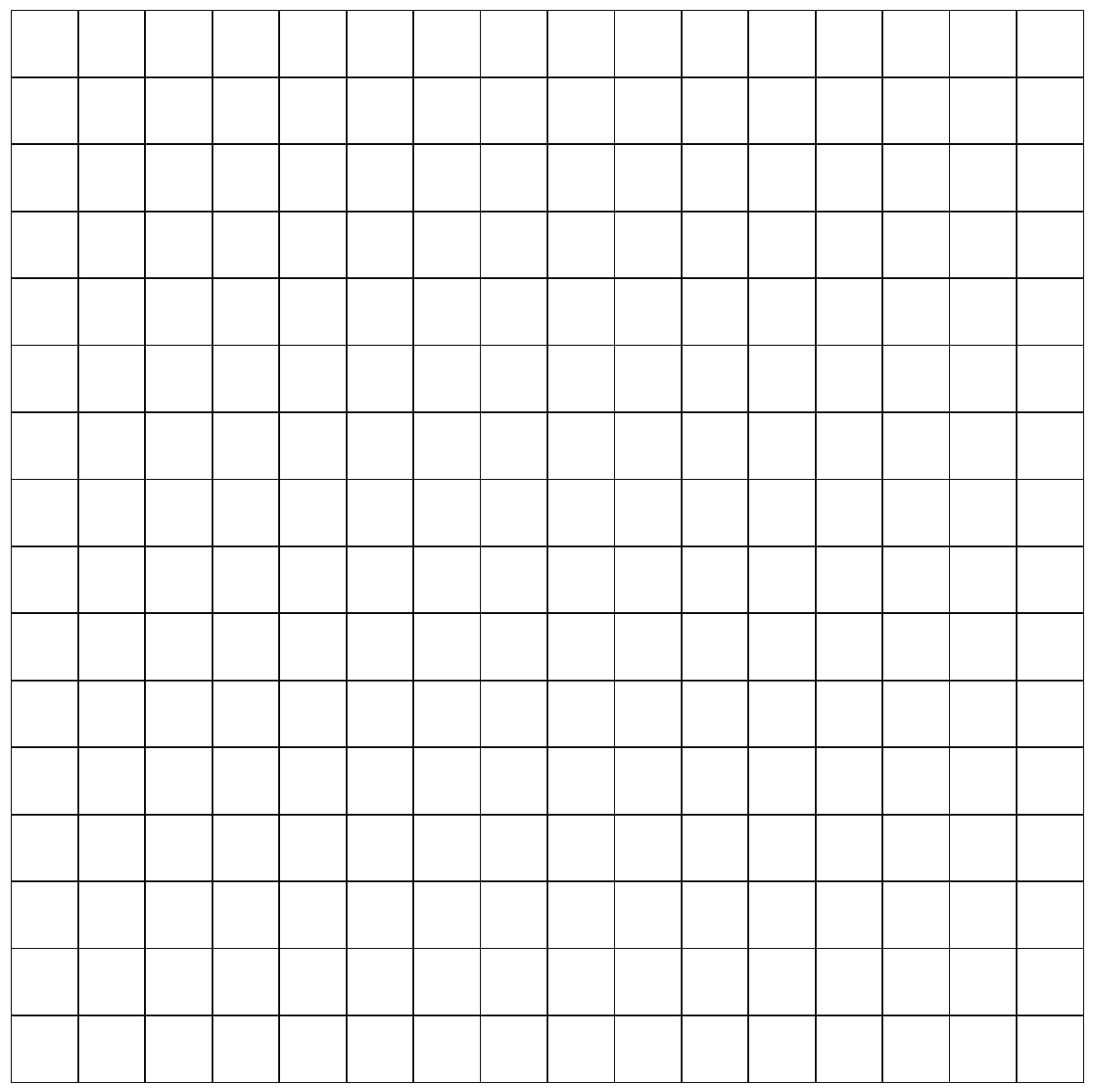}}%
\subcaptionbox{Voronoi\label{fig:voronoiMesh}}[4.7cm]{%
\includegraphics[width=4.5cm]{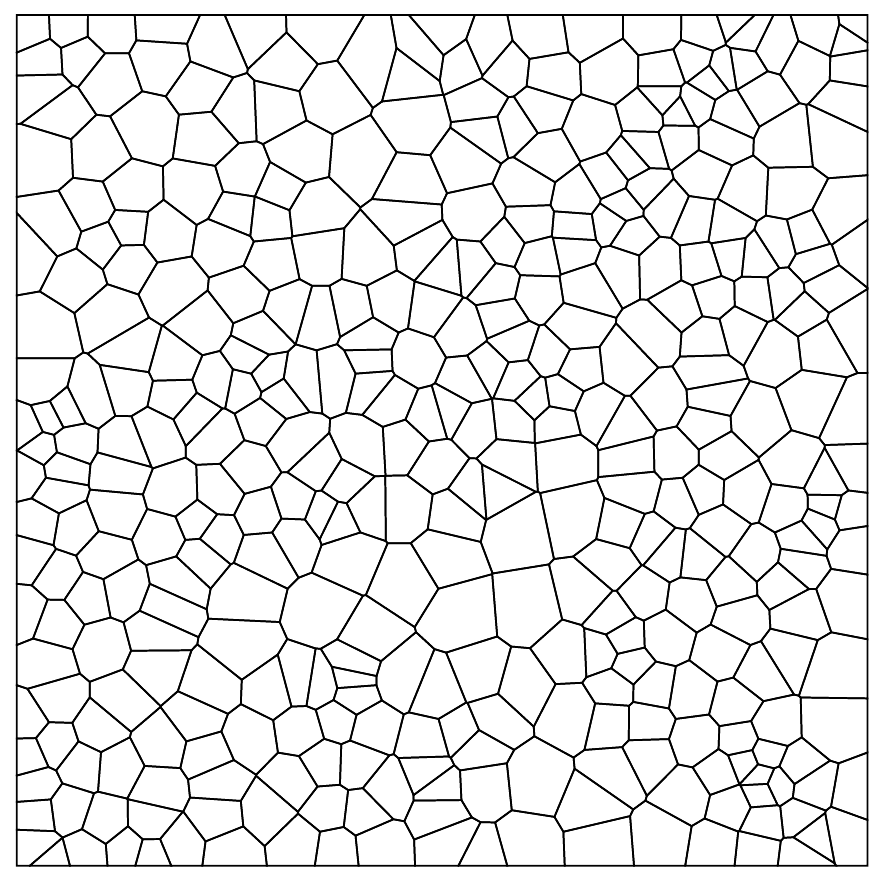}}%

\subcaptionbox{Smoothed Voronoi}[4.7cm]{%
\includegraphics[width=4.5cm]{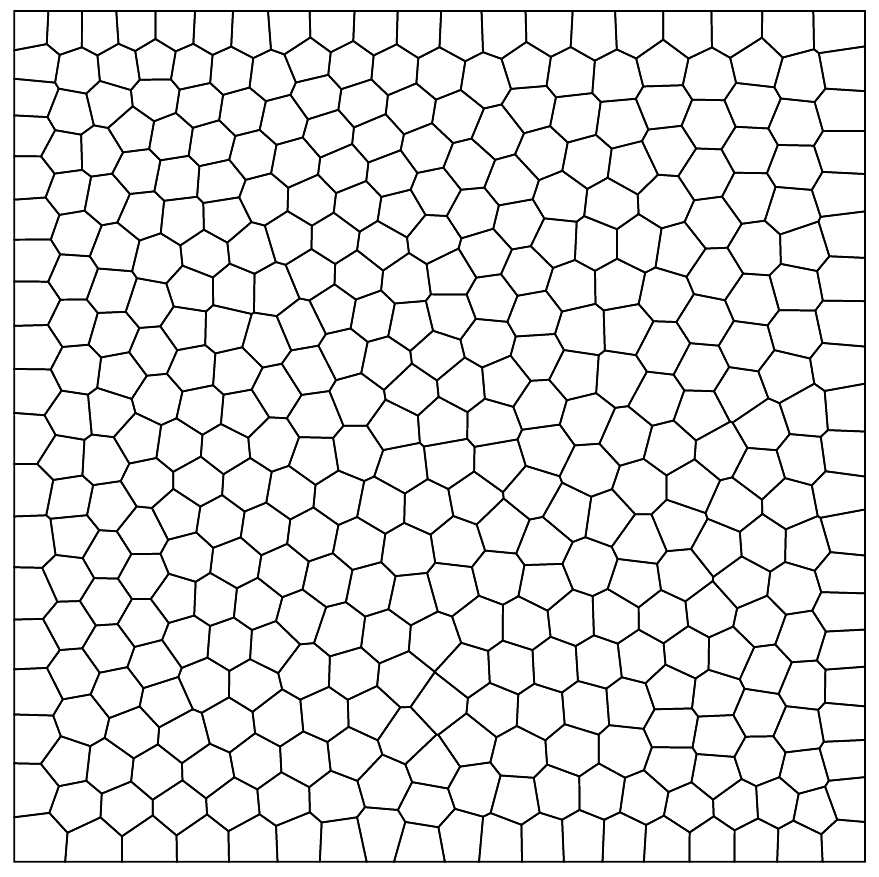}}%
\subcaptionbox{Non-convex}[4.7cm]{%
\includegraphics[width=4.5cm]{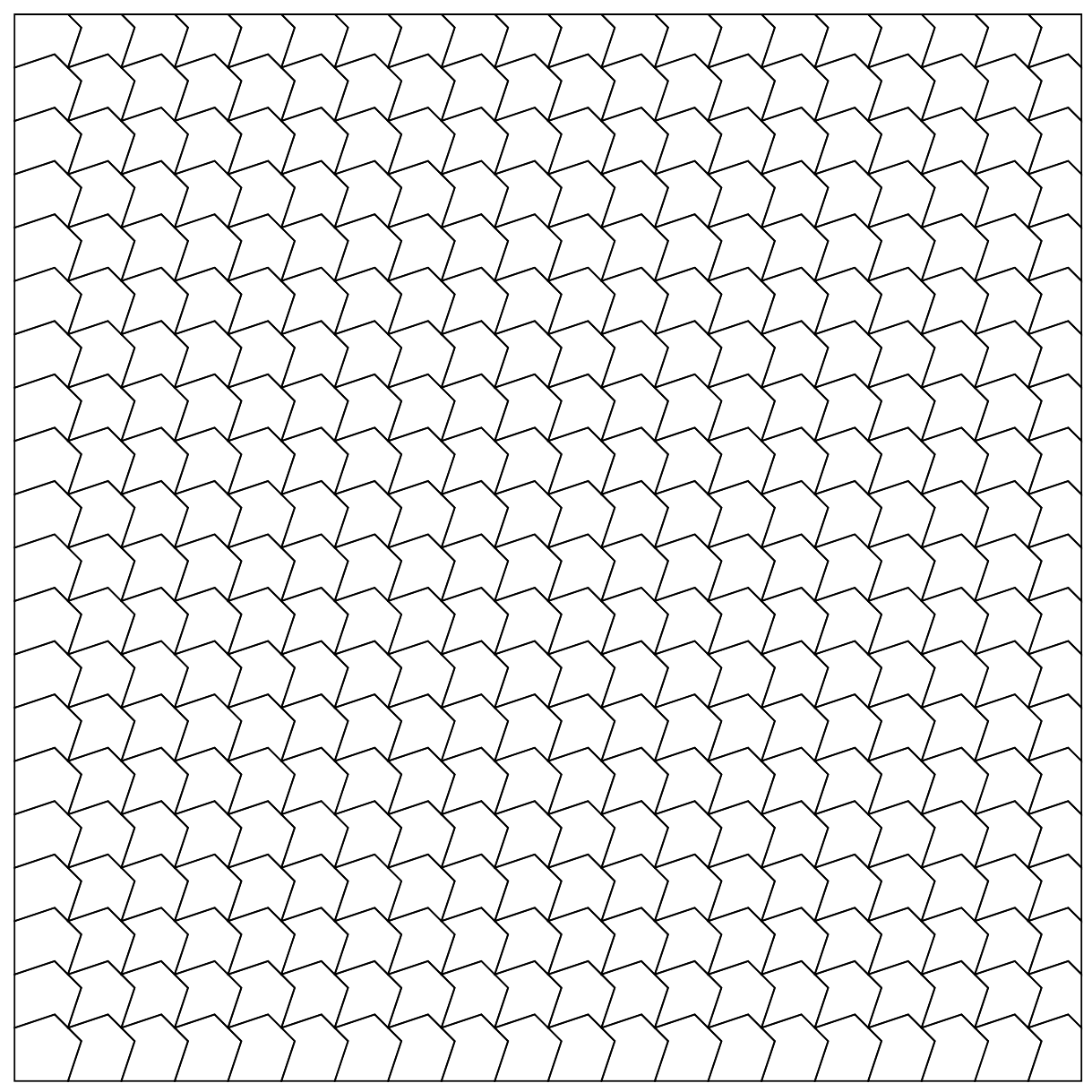}}%
\subcaptionbox{L-shaped domain\label{fig:LDomainMesh}}[4.7cm]{%
\includegraphics[width=4.5cm]{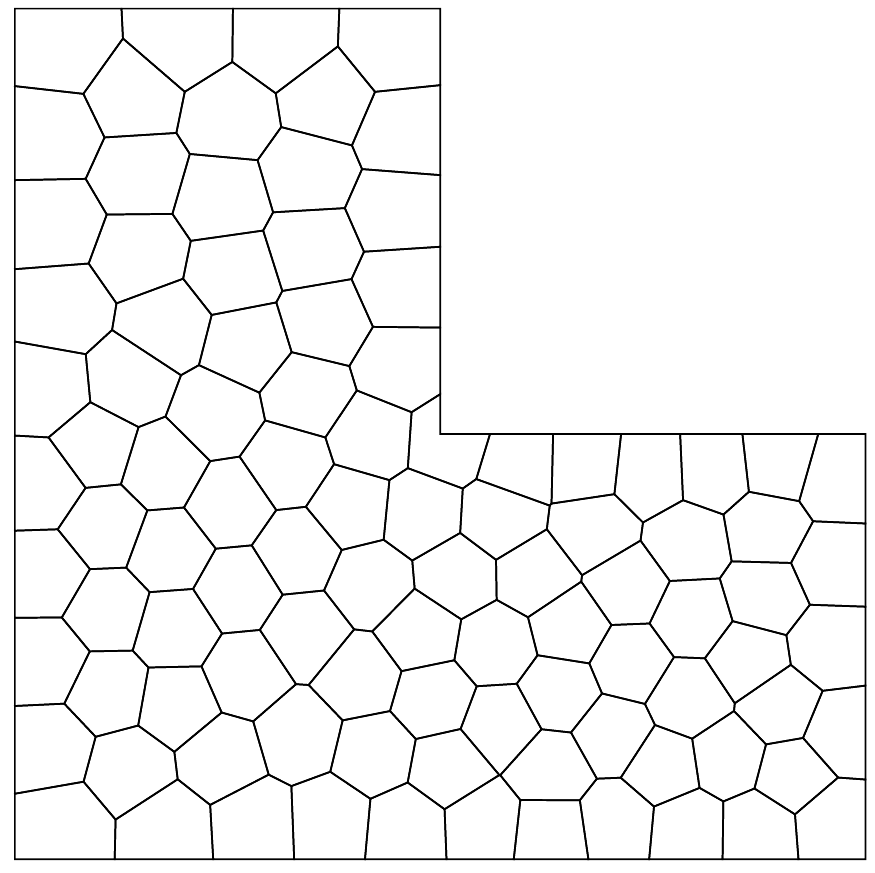}}%

\caption{The sample meshes available along with the code}
\label{fig:sampleMeshes}
\end{figure}

\subsection{The polygonal mesh}
The mesh is loaded into a structure named \texttt{mesh} from a binary MATLAB \texttt{.mat} file containing a matrix named \texttt{vertices} which specifies a mesh vertex on each row, a cell array named \texttt{elements} containing vectors of indices indicating the vertices which make up each element in an anti-clockwise order around the element, and a vector named \texttt{boundary} containing the indices of the vertices which lie on the boundary of the mesh.
Illustrated in Figure~\ref{fig:sampleMeshes} are several examples of such meshes, distributed alongside the code. 
This information can also be generated in the same format using the Voronoi mesh generator PolyMesher~\cite{Talischi:2012fx}, also written in MATLAB.

\subsection{Initialisation}
\label{sec:implementation:initialisation}
The initialisation step of the code, shown in Listing~\ref{lst:initialisation}, simply sets up various variables which will be useful to us later.
In the interests of efficiency, we use a sparse matrix \texttt{K} to represent the stiffness matrix.
In this step we also define the cell array \texttt{linear\_polynomials} containing three pairs of numbers indicating the degree of the associated polynomial in the $x$ and $y$ directions.
Thus, the index of a specific polynomial in this array is taken to be the index of the polynomial in the basis $\ME{1}$.
We note that the ordering imposed in the code coincides with the ordering in~\eqref{eq:polynomialbasis}, although this choice is arbitrary.

Some extra element-specific initialisation also takes place within the loop over all mesh element to compute various geometric properties of the element.
The vector \texttt{vert\_ids} contains the global indices of the vertices forming the element $\E$ with an anti-clockwise ordering.
As well as providing us with a means of looking up the coordinates of the vertices using the \texttt{mesh} structure, this also provides us with a very simple way to identify the global index of a particular local degree of freedom.
This is possible because the indices of the degrees of freedom of the \emph{global} virtual element space can be taken to be just the \emph{global} indices of their associated vertices.
Similarly, the \emph{local} indices of the vertices of the element $\E$ dictate the \emph{local} index of the associated degree of freedom.
Therefore, the $i^\text{th}$ entry in the vector \texttt{vert\_ids} provides the \emph{global} index of the $i^\text{th}$ \emph{local} vertex and therefore also the \emph{global} index of its associated \emph{local} degree of freedom.
Having access to this `local to global' mapping is absolutely crucial when trying to assemble the local stiffness matrix and forcing vector into their global counterparts.
When implementing more complex methods this sort of bookkeeping can quickly become very cumbersome, although here we are able to exploit the very simple arrangement of the degrees of freedom.
Because of this property, we use the variable \texttt{n\_sides} to represent both the number of sides of $\E$ and equivalently the number of local degrees of freedom of $\VhE$.

The variable \texttt{area} denotes $\abs{\E}$ and is computed using the formula
\begin{equation}
	\abs{\E} = \frac{1}{2} \left\lvert \sum_{i=1}^{\dimVhE} x_i y_{i+1} - x_{i+1} y_i \right\rvert,
	\label{eq:area}
\end{equation}
where $(x_i,y_i)$ are the coordinates of the vertex $\nu_i$, and the indexing is understood to wrap within the range 1 to $\dimVhE$.

The centroid $(x_{\E}, y_{\E})$ of the element is stored in the vector \texttt{centroid} and calculated using the usual formula:
\begin{align}
	\begin{split}
		x_{\E} &= \frac{1}{6|\E|} \sum_{i=1}^{\dimVhE} (x_i + x_{i+1}) (x_i y_{i+1} - x_{i+1} y_{i}), \\
		y_{\E} &= \frac{1}{6|\E|} \sum_{i=1}^{\dimVhE} (y_i + y_{i+1}) (x_i y_{i+1} - x_{i+1} y_{i}),
	\end{split}
	\label{eq:centroid}
\end{align}
where the indices are again to be understood to wrap within the range 1 to $\dimVhE$.
In the code, we are able to combine some of the calculations which are necessary to find the area and centroid by storing the terms of the sum~\eqref{eq:area} in the vector \texttt{area\_components}.

\begin{lstlisting}[language=matlab,firstline=8,lastline=29,firstnumber=last,caption={The initialisation process (extract from Listing~\ref{lst:vem.m})},label=lst:initialisation]
function u = vem(mesh_filename)
% Computes the virtual element solution of the Poisson problem on the specified mesh.
	mesh = load(mesh_filename); % Load the mesh from a .mat file
	n_dofs = size(mesh.vertices, 1); n_polys = 3; % Method has 1 degree of freedom per vertex
	K = sparse(n_dofs, n_dofs); % Stiffness matrix
	F = zeros(n_dofs, 1); % Forcing vector
	u = zeros(n_dofs, 1); % Degrees of freedom of the virtual element solution
	linear_polynomials = {[0,0], [1,0], [0,1]}; % Impose an ordering on the linear polynomials
	mod_wrap = @(x, a) mod(x-1, a) + 1; % A utility function for wrapping around a vector
	for el_id = 1:length(mesh.elements)
		vert_ids = mesh.elements{el_id}; % Global IDs of the vertices of this element
		verts = mesh.vertices(vert_ids, :); % Coordinates of the vertices of this element
		n_sides = length(vert_ids); % Start computing the geometric information
		area_components = verts(:,1) .* verts([2:end,1],2) - verts([2:end,1],1) .* verts(:,2);
		area = 0.5 * abs(sum(area_components));
		centroid = sum((verts + verts([2:end,1],:)) .* repmat(area_components,1,2)) / (6*area);
		diameter = 0; % Compute the diameter by looking at every pair of vertices
		for i = 1:(n_sides-1)
			for j = (i+1):n_sides
				diameter = max(diameter, norm(verts(i, :) - verts(j, :)));
			end
		end
\end{lstlisting}

\subsection{The Ritz projection and local stiffness matrix}

We focus initially on computing $\Pn \biE$ for a single basis function $\biE$.
Since $\Pn \biE \in \PE{1} \subset \VhE$, we have two different possible expansions for the projection $\Pn \biE$, either in the basis of $\VhE$ or in that of $\PE{1}$:
\begin{equation}
	\Pn \biE = \sum_{\alpha = 1}^{\dimPE} a_{i, \alpha} \ma = \sum_{j=1}^{\dimVhE} s_{i,j} \bjE.
	\label{eq:projectionExpansion}
\end{equation}
Recalling~\eqref{eq:computingPiNabla} and using the fact that $\nabla \mb$ is a constant vector for the linear polynomial $\mb$, we find that
\begin{align*}
	\sum_{\alpha=1}^{\dimPE} a_{i, \alpha} (\nabla \ma, \nabla \mb)_{0,\E} &= \sum_{j = 1}^{\dimVhE} \int_{\edge_j} \vh \n_{\edge_j} \cdot \nabla \mb \ds \\
			&= \sum_{j = 1}^{\dimVhE} \frac{\abs{\edge_j}}{2} (\biE(\nu_{j}) + \biE(\nu_{j+1})) \n_{\edge_j} \cdot \nabla \mb \\
			&= \frac{\abs{\edge_{i}}\n_{\edge_i} + \abs{\edge_{i-1}}\n_{\edge_{i-1}}}{2} \cdot \nabla \mb,
\end{align*}
for any $\mb \in \ME{1}$, where in the last equality we have used the Lagrangian property of the basis functions $\biE$ at the vertices of $\E$ to determine that $\biE$ is only non-zero on the two edges $\edge_i$ and $\edge_{i-1}$ which meet at the vertex $\nu_i$.
As shown in~\cite{Cangiani:2015ia}, this can be further simplified to
\begin{equation*}
	\sum_{\alpha=1}^{\dimPE} a_{i, \alpha} (\nabla \ma, \nabla \mb)_{0,\E} = \frac{1}{2}\abs{\widehat{\edge}_{i}}\n_{\widehat{\edge}_{i}} \cdot \nabla \mb,
\end{equation*}
where we have denoted by $\widehat{\edge}_i$ the line segment connecting the vertices $\nu_{i-1}$ and $\nu_{i+1}$, and $\n_{\widehat{\edge}_i}$ is the unit normal to $\widehat{\edge}_i$ such that $\n_{\widehat{\edge}_i} \cdot \n_{\edge_j} \geq 0$ for $j = i, i-1$.

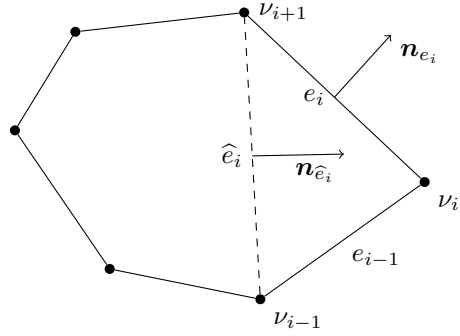
\begin{figure}
  \centering
  \begin{tikzpicture}[rotate=-20, scale=1]
    	\node[shape=coordinate] (barycentre) at (1,3.3) {};
		\node[circle, fill, scale=0.4] (v1) at (0.5,2.2) {};
		\node[circle, fill, scale=0.4, label=right:{$\nu_{i-1}$}, rotate=-45] (v3) at (2.5,2.5) {};
		\node[circle, fill, scale=0.4, label=right:{$\nu_{i}$}, rotate=-45] (v4) at (4,4.7) {};
		\node[circle, fill, scale=0.4, label=right:{$\nu_{i+1}$}, rotate=0] (v5) at (1,6) {};
		\node[circle, fill, scale=0.4] (v6) at (-1,5) {};
		\node[circle, fill, scale=0.4] (v7) at (-1.3,3.5) {};
		\node[shape=coordinate] (o1) at (0.5,0) {};
		\node[shape=coordinate] (o3) at (3.4,1.9) {};
		\node[shape=coordinate] (o4) at (4.4,5.2) {};
		\node[shape=coordinate] (o5) at (0.9,6.4) {};
		\node[shape=coordinate] (o6) at (-1.4,4.2) {};
		\node[shape=coordinate] (o7) at (-1.7,1.8) {};
		%
		\draw (v1) -- (v3) -- (v4) -- (v5) -- (v6) -- (v7) -- (v1);
		%
		%
		\draw[dashed] (v3) -- (v5);
		\node[coordinate, label=-30:{$\edge_{i-1}$}, rotate=-45] at ($0.5*(v3)+0.5*(v4)$) {};
		\node[coordinate, label=left:{$\edge_{i}$}] (eimidpoint) at ($0.5*(v5)+0.5*(v4)$) {};
		\node[coordinate, label=left:{$\widehat{\edge}_{i}$}] (fictitiousmidpoint) at ($0.5*(v5)+0.5*(v3)$) {};
		\node[shape=coordinate, label=185:{$\n_{\widehat{\edge}_{i}}$}] (innerarrowend) at ($0.5*(fictitiousmidpoint)+0.5*(v4) + (0,0.21)$) {};
		\node[shape=coordinate, label=-45:{$\n_{\edge_{i}}$}] (outerarrowend) at ($0.55*(o5)+0.55*(o4)$) {};
		\draw[->] (fictitiousmidpoint) -- (innerarrowend) {};
		\draw[->] (eimidpoint) -- (outerarrowend) {};
	\end{tikzpicture}
	\caption{An illustration of the labelling of the various geometric attributes on each element. Vertices are labelled as $\nu_j$, the edge connecting $\nu_j$ to $\nu_{j+1}$ is denoted by $\edge_j$, and $\widehat{\edge}_{i}$ denotes the lines segment connecting $\nu_{j-1}$ and $\nu_{j+1}$. The outward unit normal in each case is denoted by $\n$ with the appropriate subscript.}
	\label{fig:num:polygonRefinement}
\end{figure} 

In view of this, we introduce the matrix $\projlhsmat \in \Re^{\dimPE \times \dimPE}$ and the vector $\projrhs_i \in \Re^{\dimPE}$ such that
\begin{equation*}
	\projlhsmat_{\beta,\alpha} = (\nabla \ma, \nabla \mb)_{0,\E},
	\qquad\qquad
	\projrhs_{\beta,i} = \frac{1}{2}\abs{\widehat{\edge}_{i}}\n_{\widehat{\edge}_{i}} \cdot \nabla \mb,
\end{equation*}
to encode the conditions above.
The problem here is that the first row (and column) of $\projlhsmat$ and $\projrhs_i$ are zero, since the gradient of a constant function is zero, and therefore $\projlhsmat$ is rank deficient.
This is overcome by using the second condition in the definition~\eqref{eq:ritzdefn} of $\Pn$, defining
\begin{align}
	\widetilde{\projlhsmat}_{\beta,\alpha} = 
		\begin{cases} 
			\frac{1}{\dimVhE} \sum_{j=1}^{\dimVhE} \ma(\nu_j) &\text{ if } \beta = 1, \\
			\projlhsmat_{\beta,\alpha} &\text{ otherwise},
		\end{cases} \qquad\qquad
	\widetilde{\projrhs}_{\beta,i} = 
		\begin{cases}
			\frac{1}{\dimVhE} &\text{ if } \beta = 1, \\
			\projrhs_{\beta,i} &\text{ otherwise}.
		\end{cases}
		\label{eq:Bdefn}
\end{align}
Therefore, the coefficients $a_{i,\alpha}$ can be calculated by solving the (full rank) matrix equation
\begin{equation*}
	\widetilde{\projlhsmat} a_i = \widetilde{\projrhs}_i.
\end{equation*}

Defining the matrix $\widetilde{\projrhs} \in \Re^{\dimPE \times \dimVhE}$ such that its $i^{\text{th}}$ column is the vector $\widetilde{\projrhs}_i \in \Re^{\dimPE}$, we obtain the matrix equation
\begin{equation*}
	\widetilde{\projlhsmat} \projmat = \widetilde{\projrhs},
\end{equation*}
where $\projmat \in \Re^{\dimPE \times \dimVhE}$ is the matrix representation of the Ritz projector $\Pn$, taking a vector of coefficients of a function expressed in terms of the basis of $\VhE$ to a vector of coefficients of the basis of $\PE{1}$, and has $a_i$ as its $i^{\text{th}}$ column.

We also introduce the matrix $\Dmat \in \Re^{\dimVhE \times \dimPE}$ with $\Dmat_{i,\alpha} := \dof_i(\ma)$ as a one-way `change of basis' matrix for re-expressing polynomials in terms of the basis of $\VhE$.
Looking again at equation~\eqref{eq:computingPiNabla}, it is apparent that we can use $\Dmat$ as a shortcut to compute $\projlhsmat$ and  $\widetilde{\projlhsmat}$, since
\begin{equation*}
	\projlhsmat = \projrhs \Dmat, \qquad \widetilde{\projlhsmat} = \widetilde{\projrhs} \Dmat,
\end{equation*}
meaning we can instead compute
\begin{equation}
	\projmat = (\widetilde{\projrhs} \Dmat)^{-1} \widetilde{\projrhs}.
	\label{eq:projmatrixequation}
\end{equation}

Finally, we can compute the local stiffness matrix as
\begin{equation}
	\stiffmatE = \projmat^{\top} \projlhsmat \, \projmat + (\Id - \Dmat \projmat)^{\top} (\Id - \Dmat \projmat).
	\label{eq:localstiffnessmatrixequation}
\end{equation}
where $\Id \in \Re^{\dimVhE \times \dimVhE}$ denotes the identity matrix.
The first term of this sum corresponds to the consistency term of the discrete bilinear form, and the second term corresponds to the stabilising term.

\begin{lstlisting}[language=matlab,firstline=30,lastline=48,firstnumber=last,name=vemextracts,caption={Assembling the local stiffness matrix (extract from Listing~\ref{lst:vem.m})},label=lst:localstiffness]
		D = zeros(n_sides, n_polys); D(:, 1) = 1;
		B = zeros(n_polys, n_sides); B(1, :) = 1/n_sides;
		for vertex_id = 1:n_sides
			vert = verts(vertex_id, :); % This vertex and its neighbours
			prev = verts(mod_wrap(vertex_id - 1, n_sides), :);
			next = verts(mod_wrap(vertex_id + 1, n_sides), :);
			vertex_normal = [next(2) - prev(2), prev(1) - next(1)]; % Average of normals on edges
			for poly_id = 2:n_polys % Only need to loop over non-constant polynomials
				poly_degree = linear_polynomials{poly_id};
				monomial_grad = poly_degree / diameter; % Gradient of a linear polynomial is constant
				D(vertex_id, poly_id) = dot(vert - centroid, poly_degree) / diameter;
				B(poly_id, vertex_id) = 0.5 * dot(monomial_grad, vertex_normal);
			end
		end
		projector = (B*D) \ B; % Compute the local Ritz projector to polynomials
		stabilising_term = (eye(n_sides) - D * projector)' * (eye(n_sides) - D * projector);
		G = B*D; G(1, :) = 0;
		local_stiffness = projector' * G * projector + stabilising_term;
		K(vert_ids,vert_ids) = K(vert_ids,vert_ids) + local_stiffness; % Copy local to global
\end{lstlisting}

The code which computes the matrix form of the local Ritz projector and the local contribution of the stiffness matrix is given in Listing~\ref{lst:localstiffness}.
The first task (lines 23--24 of Listing~\ref{lst:localstiffness}) is to initialise the two matrices \texttt{D} and \texttt{B} representing their namesakes $\Dmat \in \Re^{\dimVhE \times \dimPE}$ and $\projrhs \in \Re^{\dimPE \times \dimVhE}$.
For each of these matrices, we can immediately calculate the elements corresponding to the constant polynomial basis function $\m_1$.
In the case of $\Dmat$, every element of the first column contains the value 1 since the constant function is 1 everywhere, while from~\eqref{eq:Bdefn} it may be observed that every element in the first row of $\projrhs$ is equal to $\dimVhE^{-1}$.

However, the remaining elements of $\Dmat$ and $\projrhs$ must be computed separately for each of the basis polynomials with total degree equal to 1, and for each local degree of freedom.
Computing the entries of $\Dmat$ is a straightforward task, since it just involves evaluating the basis monomials at the vertices of the polygon.

For the entries of $\projrhs$, however, we must evaluate the second expression in~\eqref{eq:Bdefn}.
The quantity $\abs{\widehat{\edge}_i} \n_{\widehat{\edge}_i}$ is simple to calculate, since it is just the vector
\begin{equation*}
	\abs{\widehat{\edge}_i} \n_{\widehat{\edge}_i} = (y_{i+1} - y_{i-1}, x_{i-1} - x_{i+1})^{\top}
\end{equation*}
due to the anti-clockwise orientation of the vertices around $\E$.
In the code, this result is stored in the variable \texttt{vertex\_normal}, so-called because it can be interpreted as a weighted normal vector at the vertex $\nu_i$.
Again, the indices here are understood to wrap within the range 1 to $\dimVhE$, and in the code this is accomplished using the utility function \texttt{mod\_wrap} which modifies the standard \texttt{mod} function to produce output in the range 1 to $\dimVhE$ rather than the range 0 to $\dimVhE-1$.
This modification to the \texttt{mod} function is necessary because arrays in MATLAB start at index 1, not 0.

To compute the entries of $\projrhs$ we also need to evaluate $\nabla \mb$.
Since $\mb$ is a linear polynomial, its gradient is simply a constant vector, and from the definition of $\ME{1}$ it is clear that
\begin{equation*}
	\nabla m_2 = (h_{\E}^{-1}, 0)^{\top}, \qquad \nabla m_3 = (0, h_{\E}^{-1})^{\top},
\end{equation*}
and hence by representing the polynomial degree of $\mb$ in the $x$ and $y$ directions using a vector with one entry 1 and one entry 0, as in the cell array \texttt{linear\_polynomials}, the gradient can be very simply calculated by just dividing by $h_{\E}$.

With the matrices $\Dmat$ and $\projrhs$ computed, we are in a position to use~\eqref{eq:projmatrixequation} to calculate the matrix $\projmat$ representing the projector $\Pn$, as shown on line 37 of Listing~\ref{lst:localstiffness}.
Consequently, with $\Dmat$, $\projrhs$ and $\projmat$ at our disposal, the local stiffness matrix can be computed as in~\eqref{eq:localstiffnessmatrixequation}.

The final step of this section is to add the elements of the local stiffness matrix to the positions in the global stiffness matrix for the corresponding global degrees of freedom.
This is accomplished on line 41 through the local to global mapping discussed in Section~\ref{sec:implementation:initialisation}.

\subsection{The local forcing vector}
To calculate the local forcing vector given in~\eqref{eq:localMatrices}, we must first compute the projection $\Po{0}\force$.
By definition, this satisfies
\begin{equation*}
	\intE \Po{0} \force \dx = \intE \force \dx,
\end{equation*}
which, because we are projecting to constants, can be simplified to
\begin{equation*}
	\Po{0}\force = \frac{1}{\abs{\E}} \intE \force \dx \approx \force(x_{\E}, y_{\E}),
\end{equation*}
where in the last relation we have used the barycentric quadrature rule on the polygon to approximate the integral.
Since we are only considering the linear virtual element method, this is sufficiently accurate to ensure the optimal order of convergence in the $H^1(\D)$ norm.
It is the use of this quadrature which produces the requirement in Section~\ref{sec:method} that the element must contain its own centroid. Clearly more general integration methods are possible (by triangulating the element, for example, or using a more advanced method such as~\cite{Sudhakar:2014dz,Sommariva:2007jl,2011CompM..47..535M,Mousavi:2010iy}, although for the sake of simplicity this is not something we pursue here.
Since each basis function of $\VhE$ is defined to be 1 at a single vertex and 0 at the others, we can express
\begin{equation*}
	\fvecE_j = \intE \overline{\bjE} \Po{0} \force \dx \approx \intE \frac{\force(x_{\E}, y_{\E})}{\dimVhE} \dx = \frac{\abs{\E} }{\dimVhE}\force(x_{\E}, y_{\E}).
\end{equation*}
The code to compute this and store the result in the appropriate positions in the global forcing vector is given in Listing~\ref{lst:localforcing}.

\begin{lstlisting}[language=matlab,firstline=49,lastline=50,firstnumber=last,name=vemextracts,caption={Implementation of the local forcing term (extract from Listing~\ref{lst:vem.m})},label=lst:localforcing]
		F(vert_ids) = F(vert_ids) + rhs(centroid) * area / n_sides;
	end
\end{lstlisting}

\subsection{Applying the boundary conditions}
The final step involves condensing the degrees of freedom associated with the boundary of the domain out of the linear system using the boundary condition, solving the resulting matrix equation, and re-applying the boundary data to the computed solution.
This part of the procedure is exactly the same as for a standard finite element method, but for completeness we briefly review the process here.

Using the subscript $B$ to denote the indices of degrees of freedom on $\dD$, and $I$ to denote those in the interior of $\D$, the matrix problem $\stiffmat U = \fvec$ can be expressed as
\begin{equation*}
	\begin{bmatrix} \stiffmat_{II} & \stiffmat_{IB} \\ \stiffmat_{BI} & \stiffmat_{BB} \end{bmatrix} 
	\begin{bmatrix} U_I \\ U_B \end{bmatrix}
	=
	\begin{bmatrix} \fvec_I \\ \fvec_B \end{bmatrix}
\end{equation*}
where $\stiffmat_{IB} = \stiffmat_{BI}^{\top}$ by the symmetry of the bilinear form.
Therefore, since $U_B$ is known, we must find $U_I$ by solving the problem
\begin{equation*}
	\stiffmat_{II} U_I = \fvec_I - \stiffmat_{IB} U_B.
\end{equation*}
This is realised on lines 46--47 of the code (shown in Listing~\ref{lst:finalsolve}) where we only store the result of solving the matrix system to the positions of the vector \texttt{u} which correspond to internal degrees of freedom, while the values of the boundary degrees of freedom are set separately on line 48.

The final line of the code uses the auxiliary function \texttt{plot\_solution}, given in Listing~\ref{lst:plot_solution.m}, to plot the vertex values of the virtual element solution and the mesh using MATLAB's \texttt{patch} function.

\begin{lstlisting}[language=matlab,firstline=51,lastline=57,firstnumber=last,caption={Solving the matrix problem and applying the boundary conditions (extract from Listing~\ref{lst:vem.m})},label=lst:finalsolve]
	boundary_vals = boundary_condition(mesh.vertices(mesh.boundary, :));
	internal_dofs = ~ismember(1:n_dofs, mesh.boundary); % Vertices which aren't on the boundary
	F = F - K(:, mesh.boundary) * boundary_vals; % Apply the boundary condition
	u(internal_dofs) = K(internal_dofs, internal_dofs) \ F(internal_dofs); % Solve
	u(mesh.boundary) = boundary_vals; % Set the boundary values
	plot_solution(mesh, u)
end
\end{lstlisting}

\section{Sample usage}
\label{sec:usage}
\begin{figure}[h]
	\centering
	\subcaptionbox{The solution produced by the default implementations of the functions $g$ and $\force$ on the Voronoi mesh in Figure~\ref{fig:voronoiMesh}.\label{fig:sampleSolution}}[7cm]{%
	\includegraphics[width=6.5cm]{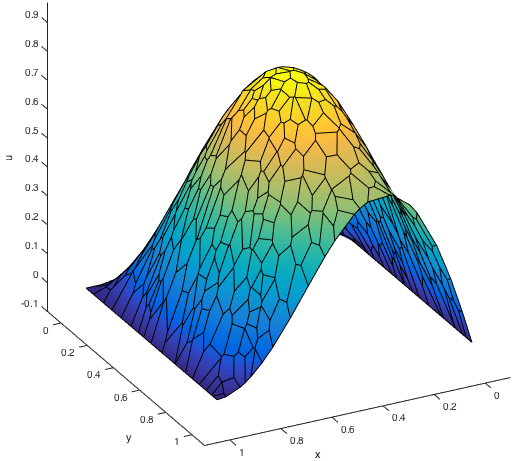}}%
	\hfill%
	\subcaptionbox{The solution to the standard example problem on the L-shaped domain mesh in Figure~\ref{fig:LDomainMesh}.\label{fig:LDomainSolution}}[7cm]{%
	\includegraphics[width=6.5cm]{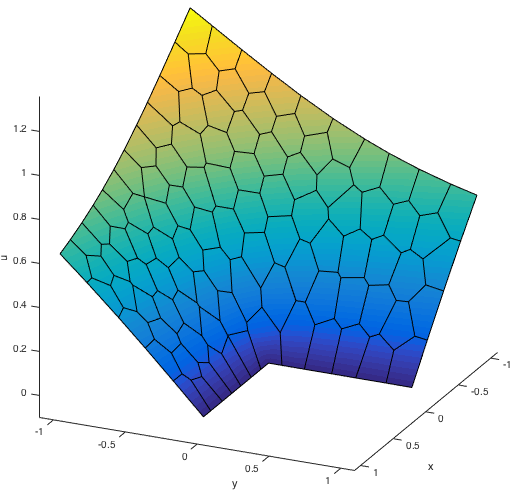}}%
\end{figure}

The default implementations of the boundary function $g$ and forcing function $\force$ are:
\begin{equation*}
	g = x y \sin(\pi x), \qquad f = 15 \sin(\pi x) \sin(\pi y),
\end{equation*}
which produces the solution shown in Figure~\ref{fig:sampleSolution}.
Computing this solution using the Voronoi mesh shown in Figure~\ref{fig:voronoiMesh} requires navigating to the directory containing the file \texttt{vem.m} within MATLAB and running
\begin{center}
	\texttt{vem('meshes/voronoi.mat');}
\end{center}
from the MATLAB Command Window.
Problems with different data can be solved simply by modifying the files \texttt{rhs.m} to change $\force$ and \texttt{boundary\_condition.m} to change $g$.
Other meshes can be used simply by specifying the path to the corresponding \texttt{.mat} file as the sole argument to the \texttt{vem} function.
For instance, setting
\begin{equation*}
	\force(r, \theta) = 0 \quad \text{ and } \quad g(r,\theta) = r^{2/3} \sin\left(\frac{2\theta - \pi}{3} \right),
\end{equation*}
where $r$ and $\theta$ are the standard polar coordinates centred at the origin, allows us to implement the standard example problem on an L-shaped domain.
Then, running \texttt{vem('meshes/L-domain.mat');} will produce the plot shown in Figure~\ref{fig:LDomainSolution}.

\section{Conclusions and extensions}
\label{sec:conclusion}
We have presented a 50-line MATLAB implementation of the linear virtual element method introduced in~\cite{VEIGA:2013wi} for solving the Poisson problem on polygonal meshes in 2 spatial dimensions, 
available to download from:
\begin{center}
\codelocation
\end{center}
alongside several example polygonal meshes.
To the best of our knowledge, this is the first publicly available implementation of the virtual element method.
It is clear from the literature surrounding the method that its capabilities extend far beyond what is presented here, although the intention behind this work is to exemplify how the method can be implemented in practice, in the simplest possible setting.
The ideas we present here can, however, be generalised to much more complicated situations by applying similar processes to compute the various required terms.
The possible extensions of this code are endless: the implementation of higher order methods, more general elliptic operators including lower order terms and non-constant coefficients~\cite{Cangiani:2015vh,BeiraodaVeiga:2014uy}, basis functions with higher global regularity properties~\cite{BeiraodaVeiga:2013ii}, mesh adaptation driven by a posteriori error indicators~\cite{Cangiani:2016ug}, or the consideration of time dependent problems~\cite{Vacca:2015jh} to name but a few.

\section*{Acknowledgements}
This work was supported by the EPSRC, which is gratefully acknowledged. The author would also like to thank H. von Wahl for correcting a mistake.

\appendix
\section{The full code}
\label{sec:fullcode}

The full code is 
available from
\begin{center}
\codelocation
\end{center}
along with the sample mesh files shown in Figure~\ref{fig:sampleMeshes}.
The software consists of four files:
\begin{itemize}
	\item \texttt{vem.m}: (see Listing~\ref{lst:vem.m}) contains the 50-line implementation of the linear virtual element method for the Poisson problem on a polygonal mesh.
	\item \texttt{rhs.m}: (see Listing~\ref{lst:rhs.m}) contains the definition of the forcing function $f$ in the model problem~\eqref{eq:modelproblem}.
	\item \texttt{boundary\_condition.m}: (see Listing~\ref{lst:boundary_condition.m}) contains the definition of the function $g$ to be used for the Dirichlet boundary condition in the model problem~\eqref{eq:modelproblem}.
	\item \texttt{plot\_solution.m}: (see Listing~\ref{lst:plot_solution.m}) a function to produce a MATLAB figure containing a plot of the approximate solution $u_h$. Note that this plot is generated using the values of $u_h$ at the vertices of the mesh, which are interpolated by the MATLAB \texttt{patch} function to produce a surface.
\end{itemize}

\clearpage
\begin{lstlisting}[language=matlab,firstline=8,caption={The main file \texttt{vem.m}},label=lst:vem.m]
function u = vem(mesh_filename)
% Computes the virtual element solution of the Poisson problem on the specified mesh.
	mesh = load(mesh_filename); % Load the mesh from a .mat file
	n_dofs = size(mesh.vertices, 1); n_polys = 3; % Method has 1 degree of freedom per vertex
	K = sparse(n_dofs, n_dofs); % Stiffness matrix
	F = zeros(n_dofs, 1); % Forcing vector
	u = zeros(n_dofs, 1); % Degrees of freedom of the virtual element solution
	linear_polynomials = {[0,0], [1,0], [0,1]}; % Impose an ordering on the linear polynomials
	mod_wrap = @(x, a) mod(x-1, a) + 1; % A utility function for wrapping around a vector
	for el_id = 1:length(mesh.elements)
		vert_ids = mesh.elements{el_id}; % Global IDs of the vertices of this element
		verts = mesh.vertices(vert_ids, :); % Coordinates of the vertices of this element
		n_sides = length(vert_ids); % Start computing the geometric information
		area_components = verts(:,1) .* verts([2:end,1],2) - verts([2:end,1],1) .* verts(:,2);
		area = 0.5 * abs(sum(area_components));
		centroid = sum((verts + verts([2:end,1],:)) .* repmat(area_components,1,2)) / (6*area);
		diameter = 0; % Compute the diameter by looking at every pair of vertices
		for i = 1:(n_sides-1)
			for j = (i+1):n_sides
				diameter = max(diameter, norm(verts(i, :) - verts(j, :)));
			end
		end
		D = zeros(n_sides, n_polys); D(:, 1) = 1;
		B = zeros(n_polys, n_sides); B(1, :) = 1/n_sides;
		for vertex_id = 1:n_sides
			vert = verts(vertex_id, :); % This vertex and its neighbours
			prev = verts(mod_wrap(vertex_id - 1, n_sides), :);
			next = verts(mod_wrap(vertex_id + 1, n_sides), :);
			vertex_normal = [next(2) - prev(2), prev(1) - next(1)]; % Average of normals on edges
			for poly_id = 2:n_polys % Only need to loop over non-constant polynomials
				poly_degree = linear_polynomials{poly_id};
				monomial_grad = poly_degree / diameter; % Gradient of a linear polynomial is constant
				D(vertex_id, poly_id) = dot(vert - centroid, poly_degree) / diameter;
				B(poly_id, vertex_id) = 0.5 * dot(monomial_grad, vertex_normal);
			end
		end
		projector = (B*D) \ B; % Compute the local Ritz projector to polynomials
		stabilising_term = (eye(n_sides) - D * projector)' * (eye(n_sides) - D * projector);
		G = B*D; G(1, :) = 0;
		local_stiffness = projector' * G * projector + stabilising_term;
		K(vert_ids,vert_ids) = K(vert_ids,vert_ids) + local_stiffness; % Copy local to global
		F(vert_ids) = F(vert_ids) + rhs(centroid) * area / n_sides;
	end
	boundary_vals = boundary_condition(mesh.vertices(mesh.boundary, :));
	internal_dofs = ~ismember(1:n_dofs, mesh.boundary); % Vertices which aren't on the boundary
	F = F - K(:, mesh.boundary) * boundary_vals; % Apply the boundary condition
	u(internal_dofs) = K(internal_dofs, internal_dofs) \ F(internal_dofs); % Solve
	u(mesh.boundary) = boundary_vals; % Set the boundary values
	plot_solution(mesh, u)
end
\end{lstlisting}
\begin{lstlisting}[language=matlab,firstline=8,caption={The file \texttt{rhs.m} which defines the forcing function $f$},label=lst:rhs.m]
function f = rhs(points)
% Evaluate the right hand side function of the PDE.
% The matrix `points` contains one point on each row.
	x = points(:, 1); y = points(:, 2);
	f = 15 * sin(pi * x) * sin(pi * y);
end
\end{lstlisting}
\begin{lstlisting}[language=matlab,firstline=8,caption={The file \texttt{boundary\_condition.m} which defines the Dirichlet boundary condition $g$},label=lst:boundary_condition.m]
function g = boundary_condition(points)
% Evaluate the boundary condition of the PDE.
% The matrix `points` contains one point on each row.
	x = points(:, 1); y = points(:, 2);
	g = (1 - x) .* y .* sin(pi * x);
end
\end{lstlisting}
\begin{lstlisting}[language=matlab,firstline=8,caption={The file \texttt{plot\_solution.m} which is used to plot the virtual element solution and polygonal mesh},label=lst:plot_solution.m]
function plot_solution(mesh, solution)
% Plot the vertex values of the virtual element solution
	figure; title('Approximate Solution');
	maxNumVertices = max(cellfun(@numel, mesh.elements));
	padFunc = @(vertList) [vertList' NaN(1,maxNumVertices-numel(vertList))];
	elements = cellfun(padFunc, mesh.elements, 'UniformOutput', false);
	elements = vertcat(elements{:});
	data = [mesh.vertices, solution];
	patch('Faces', elements, 'Vertices', data,...
							'FaceColor', 'interp', 'CData', solution / max(abs(solution)));
	axis('square')
	xlim([min(mesh.vertices(:, 1)) - 0.1, max(mesh.vertices(:, 1)) + 0.1])
	ylim([min(mesh.vertices(:, 2)) - 0.1, max(mesh.vertices(:, 2)) + 0.1])
	zlim([min(solution) - 0.1, max(solution) + 0.1])
	xlabel('x'); ylabel('y'); zlabel('u');
end
\end{lstlisting}

\bibliographystyle{acm}
\bibliography{references}

\end{document}